\documentclass[11pt]{article}
\usepackage[]{amsmath,amssymb}
\usepackage{cite}
\newtheorem{theorem}{Theorem}[section]
\newtheorem{lemma}[theorem]{Lemma}
\newtheorem{proposition}[theorem]{Proposition}
\newtheorem{corollary}[theorem]{Corollary}
\newtheorem{exAux}[theorem]{Example}

\newtheorem{Def}[theorem]{Definition}
\newenvironment{definition}{\begin{Def} \rm}{\end{Def}}
\newtheorem{Note}[theorem]{Note}
\newenvironment{note}{\begin{Note} \rm}{\end{Note}}
\newtheorem{Problem}[theorem]{Problem}

\newtheorem{Rem}[theorem]{Remark}

\newtheorem{Not}[theorem]{Notation}
\newenvironment{notation}{\begin{Not} \rm}{\end{Not}}
\newtheorem{Conj}[theorem]{Conjecture}

\newtheorem{Ass}[theorem]{Assumption}

\newenvironment{proof}{\medskip\noindent{\bf Proof.\ }}{\qed\medskip}
\newenvironment{proofof}[1]{\medskip\noindent{\bf Proof  of {#1}.\ 
}}{\qed\medskip}
\newcommand{\qed}{\hfill\mbox{$\Box$\qquad\qquad}}

\renewcommand{\b}[1]{\langle #1 \rangle}

\renewcommand{\th}{\theta}

\renewcommand{\indent}{\hspace{6mm}}

\begin{document}
\thispagestyle{empty}

\begin{center}
\LARGE \bf
\noindent
The structure of a tridiagonal pair
\end{center}

\smallskip

\begin{center}
\Large
Kazumasa Nomura and 
Paul Terwilliger\footnote{This author gratefully acknowledges 
support from the FY2007 JSPS Invitation Fellowship Program
for Reseach in Japan (Long-Term), grant L-07512.}
\end{center}

\smallskip

\begin{quote}
\small 
\begin{center}
\bf Abstract
\end{center}

\indent
Let $\mathbb{K}$ denote a field and let $V$ denote a vector space
over $\mathbb{K}$ with finite positive dimension.
We consider a pair of $\mathbb{K}$-linear transformations $A:V \to V$
and $A^*:V \to V$ that satisfy the following conditions:
(i)
each of $A,A^*$ is diagonalizable;
(ii)
there exists an ordering $\{V_i\}_{i=0}^d$ of the eigenspaces of 
$A$ such that
$A^* V_i \subseteq V_{i-1} + V_{i} + V_{i+1}$ for $0 \leq i \leq d$,
where $V_{-1}=0$ and $V_{d+1}=0$;
(iii)
there exists an ordering $\{V^*_i\}_{i=0}^\delta$ 
of the eigenspaces of $A^*$ such that
$A V^*_i \subseteq V^*_{i-1} + V^*_{i} + V^*_{i+1}$ for $0 \leq i \leq \delta$,
where $V^*_{-1}=0$ and $V^*_{\delta+1}=0$;
(iv)
there is no subspace $W$ of $V$ such that
$AW \subseteq W$, $A^* W \subseteq W$, $W \neq 0$, $W \neq V$.
We call such a pair a {\em tridiagonal pair} on $V$.
It is known that $d=\delta$ and for $0 \leq i \leq d$
the dimensions of $V_i, V_{d-i}, V^*_i, V^*_{d-i}$ coincide.
In this paper we show that the following (i)--(iv) hold
provided that $\mathbb{K}$ is algebraically closed:
(i) Each of $V_0$, $V^*_0$, $V_d$, $V^*_d$ has dimension $1$.
(ii) There exists a nondegenerate symmetric bilinear form
$\b{\;,\,}$ on $V$ such  that $\b{Au,v}=\b{u,Av}$ and $\b{A^*u,v}=\b{u,A^*v}$
for all $u,v \in V$.
(iii) There exists a unique anti-automorphism of $\text{End}(V)$ that fixes
each of $A,A^*$.
(iv) The pair $A,A^*$ is determined up to
isomorphism by the data 
$(\{\th_i\}_{i=0}^d; \{\th^*_i\}_{i=0}^d; \{\zeta_i\}_{i=0}^d)$,
where $\th_i$ (resp. $\th^*_i$) is the eigenvalue of
$A$ (resp. $A^*$) on $V_i$ (resp. $V^*_i$), and $\{\zeta_i\}_{i=0}^d$ is the
split sequence of $A,A^*$ corresponding to 
$\{\th_i\}_{i=0}^d$ and $\{\th^*_i\}_{i=0}^d$.
\end{quote}

\section{Introduction}

\indent
Throughout this paper $\mathbb{K}$  denotes a field and
$V$ denotes a vector space over $\mathbb{K}$
with finite positive dimension.

\medskip

We begin by recalling the notion of a tridiagonal pair.
We will use the following terms.
Let $\text{End}(V)$ denote the $\mathbb{K}$-algebra
consisting of all $\mathbb{K}$-linear transformations
from $V$ to $V$.
For $A \in \text{End}(V)$ and for a subspace $W \subseteq V$,
we call $W$ an {\em eigenspace} of $A$ 
whenever $W \neq 0$ and there exists $\theta \in \mathbb{K}$
such that $W=\{v \in V \,|\, Av=\theta v\}$; 
in this case $\th$ is
the {\em eigenvalue} of $A$ associated with $W$.
We say $A$ is {\em diagonalizable} whenever $V$ is spanned by
the eigenspaces of $A$.

\medskip

\begin{definition} \cite{ITT}     \label{def:TDpair}   \samepage
By a {\em tridiagonal pair} on $V$ we mean an ordered pair
$A,A^* \in \text{End}(V)$ that satisfy (i)--(iv) below:
\begin{itemize}
\item[(i)] 
Each of $A,A^*$ is diagonalizable.
\item[(ii)] 
There exists an ordering $\{V_i\}_{i=0}^d$ of the eigenspaces of 
$A$ such that
\begin{equation}              \label{eq:Astrid}
  A^* V_i \subseteq V_{i-1} + V_{i} + V_{i+1}  \qquad\qquad (0 \leq i \leq d),
\end{equation}
 where $V_{-1}=0$ and $V_{d+1}=0$.
\item[(iii)] 
There exists an ordering $\{V^*_i\}_{i=0}^\delta$ 
of the eigenspaces of $A^*$ such that
\begin{equation}           \label{eq:Atrid}
A V^*_i \subseteq V^*_{i-1} + V^*_{i} + V^*_{i+1} \qquad\qquad
            (0 \leq i \leq \delta),
\end{equation}
where $V^*_{-1}=0$ and $V^*_{\delta+1}=0$.
\item[(iv)] 
There is no subspace $W$ of $V$ such that
$AW \subseteq W$, $A^* W \subseteq W$, $W \neq 0$, $W \neq V$.
\end{itemize}
We say the tridiagonal pair $A,A^*$ is {\em over} $\mathbb{K}$.
\end{definition}

\begin{note}      \label{note:star}        \samepage
It is a common notational convention to use $A^*$ to represent the
conjugate-transpose of $A$. We are not using this convention.
In a tridiagonal pair $A,A^*$ the linear transformations $A$ and
$A^*$ are arbitrary subject to (i)--(iv) above.
\end{note}

\medskip

We refer the reader to
\cite{AC,AC2,AC3,Bas,ITT,IT:shape,IT:uqsl2hat,IT:Krawt,N:aw,N:refine,N:height1,
NT:tde,NT:sharp,NT:eddde,Vidar} for background on tridiagonal pairs.
See 
\cite{BI,BT:Borel,BT:loop,Bow,Ca,CaMT,CaW,Egge,F:RL,H:tetra,HT:tetra,
IT:non-nilpotent,IT:tetra,
IT:inverting,IT:drg,IT:loop,ITW:equitable,Leo,Koe,Mik,Mik2,
P,R:multi,R:6j,Suz,T:sub1,T:sub3,
T:qSerre,T:Kac-Moody,Z}
for related topics.

\medskip

Let $A,A^*$ denote a tridiagonal pair on $V$,
as in Definition \ref{def:TDpair}.
By \cite[Lemma 4.5]{ITT} the integers $d$ and $\delta$ from (ii), (iii)
are equal; we call this common value the {\em diameter} of the pair. 
By \cite[Corollary 5.7]{ITT}, for $0 \leq i \leq d$ the spaces $V_i$, $V^*_i$
have the same dimension; we denote
this common dimension by $\rho_i$. 
By the construction $\rho_i \neq 0$.
By \cite[Corollaries 5.7, 6.6]{ITT}
the sequence $\{\rho_i\}_{i=0}^d$ is symmetric and unimodal;
that is $\rho_i=\rho_{d-i}$ for $0 \leq i \leq d$ and
$\rho_{i-1} \leq \rho_i$ for $1 \leq i \leq d/2$.
We call the sequence $\{\rho_i\}_{i=0}^d$ the {\em shape}
of $A,A^*$.
The pair $A,A^*$ is said to be {\em sharp}
whenever $\rho_0=1$ \cite[Definition 1.5]{NT:sharp}.
By a {\em Leonard pair} we mean a tridiagonal pair with shape
$(1,1,\ldots,1)$ \cite[Definition 1.1]{T:Leonard}.
See 
\cite{Cur:mlt,Cur:spinLP,H,M:LT,NT:balanced,NT:formula,NT:det,NT:mu,
NT:span,NT:switch,NT:affine,NT:maps,T:Leonard,T:24points,T:canform,T:intro,
T:intro2,T:split,T:array,T:qRacah,T:survey,TV,V,V:AW}
for more information about Leonard pairs.

\medskip

In this paper we obtain the following four results.

\medskip

\begin{theorem}     \label{thm:sharp}   \samepage
A tridiagonal pair over an algebraically closed field is sharp.
\end{theorem}

\begin{theorem}          \label{thm:bilin}  \samepage
Let $A,A^*$ denote a sharp tridiagonal pair on $V$.
Then there exists a nonzero bilinear form $\b{\;,\,}$ on $V$ such that
$\b{Au,v}=\b{u,Av}$ and $\b{A^*u,v}=\b{u,A^*v}$ for all $u,v \in V$.
This form is unique up to multiplication by a nonzero scalar in $\mathbb{K}$.
This form is nondegenerate and symmetric.
\end{theorem}

\begin{theorem}           \label{thm:existanti}  \samepage
Let $A,A^*$ denote a sharp tridiagonal pair on $V$.
Then there exists a unique anti-automorphism $\dagger$ of $\text{\rm End}(V)$
that fixes each of $A,A^*$.
Moreover $X^{\dagger\dagger}=X$ for all $X \in \text{\rm End}(V)$.
\end{theorem}

\medskip

There is an object closely related to a tridiagonal pair
called a {\em tridiagonal system} (see Definition \ref{def:TDsystem}).
Associated with a tridiagonal system 
is a sequence of scalars called the {\em  parameter array}
(see Definition \ref{def:parray}).

\medskip

\begin{theorem}    \label{thm:determined}   \samepage
Two sharp tridiagonal systems over $\mathbb{K}$
are isomorphic if and only if they have the same parameter array.
\end{theorem}

\medskip

For Theorems \ref{thm:sharp}--\ref{thm:determined} 
the following partial results are already in the literature.
Associated with a tridiagonal pair $A,A^*$ 
is a scalar $q$ that is used to describe the
eigenvalues \cite[Lemma 8.6]{ITT}.
There is a family of tridiagonal pairs with $q=1$
said to have {\it Krawtchouk type} \cite{IT:Krawt}. 
There is a family of tridiagonal pairs with $q$ not a root of
unity said to have {\it $q$-Serre type} \cite{AC,AC2,AC3}
or {\it $q$-geometric type} \cite{IT:uqsl2hat,IT:non-nilpotent,IT:inverting}.
Theorem \ref{thm:sharp} is proven assuming $q$-geometric type 
\cite{IT:uqsl2hat}, Krawtchouk type \cite{IT:Krawt}, 
or $q$ not a root of unity \cite{IT:aug}.
Theorems \ref{thm:bilin}, \ref{thm:existanti} are proven assuming
$q$-geometric type \cite{AC3} or Krawtchouk type \cite{IT:Krawt}.
Theorem \ref{thm:determined} is proven assuming 
$q$-geometric type \cite{IT:non-nilpotent}, 
Krawtchouk type \cite{IT:Krawt},
or $q$ not a root of unity \cite{IT:aug}.
In our proof of Theorems \ref{thm:sharp}--\ref{thm:determined} 
we do not use the above partial results.
Our argument is case-free and does not refer to $q$.

\medskip

We will return to the above theorems
after a brief discussion of tridiagonal systems.

\section{Tridiagonal systems}

\indent
When working with a tridiagonal pair,
it is often convenient to consider
a closely related object called a tridiagonal system.
To define a tridiagonal system, we recall a few concepts from linear
algebra.
Let $A$ denote a diagonalizable element of $\text{End}(V)$.
Let $\{V_i\}_{i=0}^d$ denote an ordering of the eigenspaces of $A$
and let $\{\th_i\}_{i=0}^d$ denote the corresponding ordering of
the eigenvalues of $A$.
For $0 \leq i \leq d$ let $E_i:V \to V$ denote the linear transformation
such that $(E_i-I)V_i=0$ and $E_iV_j=0$ for $j \neq i$ $(0 \leq j \leq d)$.
Here $I$ denotes the identity of $\text{End}(V)$.
We call $E_i$ the {\em primitive idempotent} of $A$ corresponding to $V_i$
(or $\th_i$).
Observe that
(i) $\sum_{i=0}^d E_i = I$;
(ii) $E_iE_j=\delta_{i,j}E_i$ $(0 \leq i,j \leq d)$;
(iii) $V_i=E_iV$ $(0 \leq i \leq d)$;
(iv) $A=\sum_{i=0}^d \th_iE_i$.
Moreover
\begin{equation}         \label{eq:defEi}
  E_i=\prod_{\stackrel{0 \leq j \leq d}{j \neq i}}
          \frac{A-\th_jI}{\th_i-\th_j}.
\end{equation}
We note that each of $\{E_i\}_{i=0}^d$,
$\{A^i\}_{i=0}^d$ is a basis for the $\mathbb{K}$-subalgebra
of $\text{End}(V)$ generated by $A$.
Moreover $\prod_{i=0}^d(A-\th_iI)=0$.

Now let $A,A^*$ denote a tridiagonal pair on $V$.
An ordering of the eigenspaces of $A$ (resp. $A^*$)
is said to be {\em standard} whenever it satisfies 
\eqref{eq:Astrid} (resp. \eqref{eq:Atrid}). 
We comment on the uniqueness of the standard ordering.
Let $\{V_i\}_{i=0}^d$ denote a standard ordering of the eigenspaces of $A$.
Then the ordering $\{V_{d-i}\}_{i=0}^d$ is standard and no other ordering
is standard.
A similar result holds for the eigenspaces of $A^*$.
An ordering of the primitive idempotents of $A$ (resp. $A^*$)
is said to be {\em standard} whenever
the corresponding ordering of the eigenspaces of $A$ (resp. $A^*$)
is standard.

\medskip

\begin{definition} \cite[Definition 2.1]{ITT}  \label{def:TDsystem}  \samepage
By a {\em tridiagonal system} on $V$ we mean a sequence
\[
   \Phi= (A;\{E_i\}_{i=0}^d;A^*;\{E^*_i\}_{i=0}^d)
\]
that satisfies (i)--(iii) below.
\begin{itemize}
\item[(i)]
$A,A^*$ is a tridiagonal pair on $V$.
\item[(ii)]
$\{E_i\}_{i=0}^d$ is a standard ordering
of the primitive idempotents of $A$.
\item[(iii)]
$\{E^*_i\}_{i=0}^d$ is a standard ordering
of the primitive idempotents of $A^*$.
\end{itemize}
We say $\Phi$ is {\em over} $\mathbb{K}$.
\end{definition}

\medskip

We will use the following notation.

\medskip

\begin{notation}   \label{notation}   \samepage
let $\Phi=(A;\{E_i\}_{i=0}^d;A^*;\{E^*_i\}_{i=0}^d)$ denote
a tridiagonal system on $V$.
For $0 \leq i \leq d$ let $\th_i$ (resp. $\th^*_i$) denote the
eigenvalue of $A$ (resp. $A^*$) associated with the eigenspace
$E_iV$ (resp. $E^*_iV$).
We observe that $\{\th_i\}_{i=0}^d$
(resp. $\{\th^*_i\}_{i=0}^d$) are mutually distinct
and contained in $\mathbb{K}$.
$\Phi$ is said to be {\em sharp} whenever
the tridiagonal pair $A,A^*$ is sharp.
\end{notation}

\medskip

To prove Theorems \ref{thm:sharp}--\ref{thm:determined} 
we will introduce a certain algebra $T$ defined by generators and relations.
The algebra $T$ is reminiscent of the algebra ${\cal T}$
introduced by E. Egge \cite[Definition 4.1]{Egge}.
Our definition of $T$ is motivated by the following observation.

\medskip

\begin{lemma} {\rm \cite[Lemma 2.5]{NT:eddde}}  \label{lem:trid}   \samepage
With reference to Notation {\rm \ref{notation}}
the following {\rm (i)}, {\rm (ii)} hold for $0 \leq i,j,k \leq d$.
\begin{itemize}
\item[\rm (i)]
$E^*_iA^kE^*_j=0\;$ if $\;k<|i-j|$.
\item[\rm (ii)]
$E_i{A^*}^kE_j=0\;$ if $\;k<|i-j|$.
\end{itemize}
\end{lemma}

\begin{definition}           \label{def:T}       \samepage
With reference to Notation \ref{notation}, let $T$ denote
the unital associative $\mathbb{K}$-algebra defined by 
generators $(a;\{e_i\}_{i=0}^d;a^*;\{e^*_i\}_{i=0}^d)$ and relations
\begin{equation}                            \label{eq:eiej}
  e_ie_j=\delta_{i,j}e_i, \qquad 
  e^*_ie^*_j=\delta_{i,j}e^*_i \qquad\qquad 
  (0 \leq i,j \leq d),
\end{equation}
\begin{equation}                            \label{eq:sumei}
  \sum_{i=0}^d e_i=1, \qquad\qquad
  \sum_{i=0}^d e^*_i=1,
\end{equation}
\begin{equation}                                  \label{eq:sumthiei}
   a = \sum_{i=0}^d \th_ie_i, \qquad\qquad
   a^* = \sum_{i=0}^d \th^*_i e^*_i,
\end{equation}
\begin{equation}                                     \label{eq:esiakesj}
 e^*_i a^k e^*_j = 0  \qquad \text{if $\;k<|i-j|$}
                  \qquad\qquad (0 \leq i,j,k \leq d),
\end{equation}
\begin{equation}                                      \label{eq:eiaskej}
 e_i {a^*}^k e_j = 0  \qquad \text{if $\;k<|i-j|$} 
                  \qquad\qquad (0 \leq i,j,k \leq d). 
\end{equation}
Let $D$ (resp. $D^*$) denote the $\mathbb{K}$-subalgebra 
of $T$ generated by $a$ (resp. $a^*$).
\end{definition}

\medskip

\begin{lemma}   \label{lem:Tmodule}   \samepage
With reference to Notation {\rm \ref{notation}} and
Definition {\rm \ref{def:T}},
there exists a unique $T$-module structure on $V$ such that
$a$, $\{e_i\}_{i=0}^d$, $a^*$, $\{e^*_i\}_{i=0}^d$ act on $V$
as $A$, $\{E_i\}_{i=0}^d$, $A^*$, $\{E^*_i\}_{i=0}^d$ respectively.
This $T$-module structure is irreducible.
\end{lemma}

\begin{proof}
The $T$-module structure exists by Lemma \ref{lem:trid} and
the equations above line \eqref{eq:defEi}.
The $T$-module structure is irreducible by
Definition \ref{def:TDpair}(iv).
\end{proof}

\medskip

As we will see, Theorems \ref{thm:sharp}--\ref{thm:determined}
are all implied by the following result.

\medskip

\begin{theorem}   \label{thm:es0Tes0}    \samepage
With reference to Definition {\rm \ref{def:T}}
the following {\rm (i)}--{\rm (iii)} hold.
\begin{itemize}
\item[\rm (i)]
The algebra $e^*_0 T e^*_0$ is generated by $e^*_0 D e^*_0$.
\item[\rm (ii)]
The elements of $e^*_0 D e^*_0$ mutually commute.
\item[\rm (iii)]
The algebra $e^*_0 T e^*_0$ is commutative.
\end{itemize}
\end{theorem}

\medskip

In this paper we will prove Theorem \ref{thm:es0Tes0} and use it
to prove Theorems \ref{thm:sharp}--\ref{thm:determined}.
Our proof of Theorem \ref{thm:es0Tes0} is summarized as follows. 
For any subsets $X$, $Y$ of $T$ let $XY$ denote the 
$\mathbb{K}$-subspace of $T$ spanned by $\{xy \,|\, x \in X,\, y \in Y\}$.
We first show that for $0 \leq r,s\leq d$,
\begin{equation}            \label{eq:esrDDsDessaux}
   e^*_r D D^* D e^*_s 
  = \sum_{t=0}^{\lfloor (r+s)/2 \rfloor} 
         e^*_r D e^*_t D e^*_s,
\end{equation}
where $\lfloor m \rfloor$ denotes the greatest integer less 
than or equal to $m$. 
Using \eqref{eq:esrDDsDessaux} we show that for all integers $n \geq 0$,
\[
  e^*_0 D D^* D D^* D \cdots D e^*_0  \qquad\qquad  \text{($n$ $D$'s)}
\]
is equal to 
\[
  e^*_0 D e^*_0 D e^*_0 D \cdots D e^*_0  \qquad\qquad \text{($n$ $D$'s)}.
\]
Theorem \ref{thm:es0Tes0}(i) is a routine consequence of this.
To get Theorem \ref{thm:es0Tes0}(ii) we will invoke
\cite[Theorem 8.3]{NT:eddde}.

\section{The algebra $T$}

\indent
In this section we establish some  basic properties of $T$.

\medskip

\begin{definition}    \label{def:p}  \samepage
With reference to Notation {\rm \ref{notation}}
and Definition {\rm \ref{def:T}}, 
let $p: T \to \text{End}(V)$ denote the $\mathbb{K}$-algebra
homomorphism induced by the $T$-module structure in Lemma \ref{lem:Tmodule}.
Note that $p(a)=A$, $p(a^*)=A^*$ and
$p(e_i)=E_i$, $p(e^*_i)=E^*_i$ for $0 \leq i \leq d$. 
\end{definition}

\begin{lemma}  \label{lem:injection}  \samepage
With reference to Notation {\rm \ref{notation}}
and Definition {\rm \ref{def:p}}
the following {\rm (i)}, {\rm (ii)} hold.
\begin{itemize}
\item[\rm (i)]
The restriction of $p$ to $D$ is an injection.
Moreover each of $\{a^i\}_{i=0}^d$, $\{e_i\}_{i=0}^d$ is a basis for $D$.
\item[\rm (ii)]
The restriction of $p$ to $D^*$ is an injection.
Moreover each of $\{{a^*}^i\}_{i=0}^d$, $\{e^*_i\}_{i=0}^d$ is a basis for $D^*$.
\end{itemize}
\end{lemma}

\begin{proof}
(i):
Let $D'$ denote the subspace of $T$ spanned by $\{e_i\}_{i=0}^d$.
The space $D'$ is a subalgebra of $T$ by 
\eqref{eq:eiej}, \eqref{eq:sumei}.
The space $D'$ contains $D$ by \eqref{eq:sumthiei}, 
so $D$ has dimension at most $d+1$.
The elements $\{a^i\}_{i=0}^d$ are linearly independent
since $\{A^i\}_{i=0}^d$ are linearly independent;
therefore the dimension of $D$ is at least $d+1$.
By these comments the dimension of $D$ is $d+1$ so 
$D=D'$ and the result follows.

(ii): Similar to the proof of (i) above.
\end{proof}

\begin{lemma}           \samepage
With reference to Definition {\rm \ref{def:T}},
\begin{equation}
 ae_i=\th_ie_i, \qquad\qquad
 a^*e^*_i=\th^*_ie^*_i    \qquad\qquad (0 \leq i \leq d),
\end{equation}
\begin{equation}         \label{eq:eiesi}
 e_i=\prod_{\stackrel{0\leq j\leq d}{j\neq i}}
      \frac{a-\th_j}{\th_i-\th_j},             \qquad\qquad
 e^*_i=\prod_{\stackrel{0\leq j\leq d}{j\neq i}}
      \frac{a^*-\th^*_j}{\th^*_i-\th^*_j}
 \qquad\qquad (0 \leq i \leq d),
\end{equation}
\begin{equation}         \label{eq:prodathi}
 \prod_{i=0}^d (a-\th_i)=0,    \qquad\qquad
  \prod_{i=0}^d (a^*-\th^*_i)=0.
\end{equation}
\end{lemma}

\begin{proof}
Routinely verified using \eqref{eq:eiej}--\eqref{eq:sumthiei}.
\end{proof}

\medskip

We now give some bases for $D$ and $D^*$ that will be useful
later in the paper. We will state our results for $D$;
of course similar results hold for $D^*$.

\medskip

\begin{lemma}      \label{lem:replace}  \samepage
With reference to Definition {\rm \ref{def:T}} consider 
the following basis for $D$:
\begin{equation}   \label{eq:E0Ed}
  e_0,e_1,\ldots,e_d.
\end{equation}
For $0 \leq n \leq d$, if we replace any $(n+1)$-subset of \eqref{eq:E0Ed}
by $1,a,a^2, \ldots, a^n$ then the result is still a basis for $D$.
\end{lemma}

\begin{proof}
Let $\Delta$ denote a $(n+1)$-subset of $\{0,1,\ldots, d\}$
and let $\overline{\Delta}$ denote the complement of
$\Delta$ in $\{0,1,\ldots, d\}$. We show
\begin{equation}                  \label{eq:basisaux}     
   \{a^i\}_{i=0}^n \cup \{e_i\}_{i \in \overline{\Delta}}       
\end{equation}
is a basis for $D$. The number of elements in \eqref{eq:basisaux}
is $d+1$ and this equals the dimension of $D$. Therefore it
suffices to show the elements \eqref{eq:basisaux} span $D$. 
Let $S$ denote the subspace of $D$ spanned by \eqref{eq:basisaux}. 
To show $D = S$ we show $e_i \in S$ for $i \in \Delta$. 
For $0 \leq j \leq n$ we have $a^j = \sum_{i=0}^d \theta^j_i e_i$.
In these equations we rearrange the terms to find
\begin{equation}           \label{eq:basisaux2}
  \sum_{i \in \Delta} \th^j_i e_i \in  S \qquad\qquad  (0 \leq j \leq n).
\end{equation}
In the linear system \eqref{eq:basisaux2} the coefficient matrix is 
Vandermonde and hence nonsingular. Therefore $e_i \in S$ for $i \in \Delta$.
Now $S=D$ and the result follows.
\end{proof}

\section{The proof of Theorem \ref{thm:es0Tes0}}

\indent

In this section we prove Theorem 2.6. We begin with some definitions.

\medskip

\begin{definition}    \label{def:pi}   \samepage
With reference to Definition \ref{def:T} we consider the tensor product 
$D \otimes D^* \otimes D$
where $\otimes = \otimes_{\mathbb{K}}$.
We fix integers $r,s$ $(0 \leq r,s\leq d)$ and 
define a $\mathbb{K}$-linear transformation
\[
 \pi : \qquad
 \begin{array}{ccc}
    D \otimes D^* \otimes D & \qquad \to \qquad & T \\
   X \otimes Y \otimes Z  & \qquad \mapsto \qquad & e^*_rXYZe^*_s
 \end{array}
\]
We note that the image of $\pi$ is 
$e^*_r D D^* D e^*_s$.
\end{definition}

\begin{definition}        \label{def:R}      \samepage
With reference to Definitions \ref{def:T} and \ref{def:pi},
let $R$ denote the sum of the following three subspaces of
$D \otimes D^* \otimes D$:
\begin{equation}        \label{eq:L}
 \text{Span}\{a^i \otimes e^*_j \,|\, 0 \leq i,j \leq d,\, i < |r-j| \}
    \otimes D,
\end{equation}
\begin{equation}        \label{eq:R} 
 D \otimes 
 \text{Span}\{e^*_j \otimes a^i \,|\, 0 \leq i,j \leq d, i <|s-j| \}, 
\end{equation}
\begin{equation}        \label{eq:M}
 \text{Span}\{e_i \otimes a^{*t} \otimes e_j
                    \,|\, 0 \leq i,j,t \leq d, \, t<|i-j| \}.
\end{equation}
\end{definition}

\begin{lemma}   \label{lem:kernel}   \samepage
With reference to Definitions {\rm \ref{def:pi}} and {\rm \ref{def:R}},
the space $R$ is contained in the kernel of $\pi$.
\end{lemma}

\begin{proof}
Routinely obtained using \eqref{eq:esiakesj}, \eqref{eq:eiaskej}.
\end{proof}

\begin{proposition}   \label{prop:EsrDDsDEss}   \samepage
With reference to Definitions {\rm \ref{def:T}} and  {\rm \ref{def:pi}}, 
the space  $D \otimes D^* \otimes D$
is the sum of the space $R$ from Definition {\rm \ref{def:R}}
and the following space:
\begin{equation}   \label{eq:defS}
   D \otimes 
  \text{\rm Span}\{e^*_t \,|\, 0 \leq t \leq \lfloor (r+s)/2 \rfloor \}
  \otimes D.
\end{equation}
\end{proposition}

\begin{proof}
We assume $r \leq s$; for the case $r \geq s$ we can proceed
in a similar way.
We set
\[
  m = \left\lfloor \frac{r+s}{2} \right\rfloor.
\]
Note that
\begin{equation}     \label{eq:2mr+s}
     r+s-1 \leq 2m \leq r+s.
\end{equation}
Let S denote the sum of $R$ and \eqref{eq:defS}. 
We show $S = D \otimes D^* \otimes D$. 
The following notation will be useful. 
For integers $i,j$ let $\{i,\ldots, j\}$ denote the set of integers
$\{h \,|\, i \leq h \leq j\}$.
Note that $\{i,\ldots,j\}$ is empty if $i>j$.

We first claim that the following (i)--(iii) hold
for a subset $\Delta$ of $\{0,\ldots,d\}$
and for $0 \leq i,j,k \leq d$.
\begin{itemize}
\item[\rm (i)]
Assume $|\Delta|>d-|r-k|$ and
$\text{\rm Span}\{e_h \,|\, h \in \Delta\} \otimes e^*_k \otimes e_j
      \subseteq S$.
Then $D \otimes e^*_k \otimes e_j \subseteq S$.
\item[\rm (ii)]
Assume $|\Delta|>d-|i-j|$ and
$e_i \otimes \text{\rm Span}\{e^*_h \,|\, h \in \Delta\} \otimes e_j
      \subseteq S$.
Then $e_i \otimes D^* \otimes e_j \subseteq S$.
\item[\rm (iii)]
Assume $|\Delta|>d-|k-s|$ and
$e_i \otimes e^*_k \otimes \text{\rm Span}\{e_h \,|\, h \in \Delta\}
      \subseteq S$.
Then $e_i \otimes e^*_k \otimes D \subseteq S$.
\end{itemize}
To prove (i), observe by Lemma \ref{lem:replace} 
that the following set spans $D$:
\[
   \{e_h \,|\, h \in \Delta\} \cup \{1,a,a^2,\ldots,a^{|r-k|-1}\}.
\]
Therefore $D \otimes e^*_k \otimes e_j$ is the sum of the following subspaces:
\begin{align}
 \text{Span}\{e_h \,|\, h \in \Delta\} & \otimes e^*_k \otimes e_j,
                                               \label{eq:spa1}  \\
 \text{Span}\{1,a,a^2,\ldots,a^{|r-k|-1}\} & \otimes e^*_k \otimes e_j.
                                               \label{eq:spa2}
\end{align}
The space \eqref{eq:spa1} is contained in $S$ by the assumption.
The space \eqref{eq:spa2} is contained in $S$ by \eqref{eq:L}.
Thus $D \otimes e^*_k \otimes e_j$ is contained in $S$ and (i) follows.
Parts (ii), (iii) of the claim are similarly proved.

To obtain $S = D \otimes D^* \otimes D$ we show that for $m+1 \leq k \leq d$,
\begin{equation}                      \label{eq:aux1}
  D \otimes e^*_k \otimes D \subseteq S.
\end{equation}
Our proof is by induction on $k=d,d-1,\ldots,m+1$.
Let $k$ be given. By induction we assume
\begin{equation}                      \label{eq:aux2}
 D \otimes e^*_\ell \otimes D \subseteq S  \qquad\qquad  (k< \ell \leq d).
\end{equation}

We claim that
\begin{equation}                       \label{eq:aux3}
  e_i \otimes D^* \otimes e_j \subseteq S \qquad\qquad 
                                  (0 \leq i,j\leq d,\;  k-m \leq |i-j|).  
\end{equation}
To prove the claim, let $i,j$ be given.
By \eqref{eq:defS} and by \eqref{eq:aux2} we have
\[
e_i \otimes \text{Span}\{e^*_h \,|\, h \in \Delta\} \otimes e_j \subseteq S,
\]
where
\[
 \Delta=\{0,\ldots,m\} \cup \{k+1,\ldots,d\}.
\]
We have $|\Delta|=m+1+d-k$ and $k-m\leq |i-j|$ so $|\Delta|>d-|i-j|$. 
Now $e_i \otimes D^* \otimes e_j \subseteq S$ by (ii), 
and the claim follows.

Now to show \eqref{eq:aux1}, we prove that
$D \otimes e^*_k \otimes e_\ell \subseteq S$ for $0 \leq \ell \leq d$.
We will obtain this by an inductive argument that involves simultaneously
showing a second assertion. 
Specifically, we show that for $0 \leq \ell \leq d$ both
\begin{align}
 D \otimes e^*_{k} \otimes e_j &\subseteq S
    \qquad \qquad (\ell     \leq j \leq d),    \label{eq:Aell} \\
 e_i \otimes e^*_{k} \otimes D &\subseteq S
    \qquad\qquad (\ell+m-r+1 \leq i \leq d).           \label{eq:Bell}
\end{align}
Our argument is by induction on $\ell=d,d-1,\ldots,0$.
Let $\ell$ be given. By induction we may assume
\begin{align}
 D \otimes e^*_{k} \otimes e_j & \subseteq S
    \qquad \qquad (\ell+1 \leq j \leq d),    \label{eq:Aell-1} \\
 e_i \otimes e^*_{k} \otimes D & \subseteq S
    \qquad\qquad (\ell+m-r+2 \leq i \leq d).           \label{eq:Bell-1}
\end{align}
We show \eqref{eq:Bell}.
We assume $\ell+m-r+1 \leq d$; otherwise we have nothing to prove.
By \eqref{eq:Bell-1} it suffices to show
\begin{equation}    \label{eq:target1}
  e_{\ell+m-r+1} \otimes e^*_{k} \otimes D \subseteq S.
\end{equation}
To show this, we invoke (iii) with
\[
  \Delta=\{0,\ldots,\ell+2m-k-r+1\} \cup \{\ell+1,\dots,d\}.
\]
We show
\begin{equation}           \label{eq:sp3a}
 e_{\ell+m-r+1} \otimes e^*_{k} \otimes \text{Span}\{e_{\ell+1},\ldots,e_d\}
    \subseteq S,
\end{equation}
\begin{equation}            \label{eq:sp1a}
 e_{\ell+m-r+1} \otimes e^*_{k} 
      \otimes \text{Span}\{e_0,e_1,\ldots,e_{\ell+2m-k-r+1}\}  \subseteq S,
\end{equation}
\begin{equation}    \label{eq:delta1}
   |\Delta|>d-|k-s|.
\end{equation}
Line \eqref{eq:sp3a} follows from \eqref{eq:Aell-1}. 
To obtain \eqref{eq:sp1a} we fix $j$ such that
\begin{equation}          \label{eq:sp1ai}
  0 \leq j \leq \ell+2m-k-r+1,
\end{equation}
and show $e_{\ell+m-r+1} \otimes e^*_{k} \otimes e_j \in S$. 
To do this we invoke \eqref{eq:aux3} with $i=\ell+m-r+1$.
It suffices to show
\begin{equation}          \label{eq:aux1a}
 k-m \leq |\ell+m-r+1-j|.
\end{equation}
Observe $j < \ell+m-r+1$ by $m+1\leq k$ and \eqref{eq:sp1ai},
so $|\ell+m-r+1-j |=\ell+m-r+1-j$.
Using this and \eqref{eq:sp1ai} we obtain \eqref{eq:aux1a},
and  \eqref{eq:sp1a} follows.
Next we show \eqref{eq:delta1}. 
First assume
\begin{equation}    \label{eq:delta1aux1}
  \ell+1<k+r-2m.
\end{equation}
In this case $|\Delta|=d-\ell$. 
By \eqref{eq:2mr+s} we have $r+s-1\leq 2m$. 
By this and \eqref{eq:delta1aux1} we obtain $\ell < k-s$. 
By assumption $\ell \geq 0$ so $k> s$; therefore
$|k-s|=k-s$ and \eqref{eq:delta1} follows.
Next assume
\[
  0<k+r-2m \leq \ell+1. 
\]
In this case $|\Delta|= d+2m-k-r+2$. 
If $s\leq k$ then $|k-s|=k-s$ and \eqref{eq:delta1} 
holds using \eqref{eq:2mr+s}.
If $s> k$ then $|k-s|=s-k$ and this yields \eqref{eq:delta1} since
\[
     |\Delta|-d+|k-s|  =  2(s-k)+2m+2-r-s > 0
\]
in view of \eqref{eq:2mr+s}.
Next assume
\[
  k+r-2m \leq 0.
\]
In this case $|\Delta|=d+1$ so \eqref{eq:delta1} holds. In any case
\eqref{eq:delta1} holds, and this completes the proof of \eqref{eq:Bell}.

We now show \eqref{eq:Aell}.
By \eqref{eq:Aell-1} it suffices to show
\begin{equation}    \label{eq:target2}
  D \otimes e^*_k \otimes e_{\ell} \subseteq S.
\end{equation}
To show this we invoke (i) with
\[
 \Delta=\{0,\ldots,\ell-k+m\} \cup \{\ell+m-r+1,\ldots,d\}.
\]
It suffices to show
\begin{equation}                 \label{eq:sp3b} 
 \text{Span}\{e_{\ell+m-r+1},e_{\ell+m-r+2},\ldots,e_d\} \otimes e^*_{k} 
    \otimes  e_{\ell} \subseteq S,
\end{equation}
\begin{equation}                   \label{eq:sp1b}
 \text{Span}\{e_0,e_1,\ldots,e_{\ell-k+m}\} \otimes e^*_{k}  \otimes 
         e_{\ell} \subseteq S,
\end{equation}
\begin{equation}                 \label{eq:delta2}
  |\Delta|>d-|r-k|.
\end{equation}
Line \eqref{eq:sp3b} follows from \eqref{eq:Bell}.
To obtain \eqref{eq:sp1b} we fix $i$ such that
\begin{equation}     \label{eq:delta2aux1}
  0 \leq i \leq \ell-k+m
\end{equation}
and show $e_i \otimes e^*_k \otimes e_{\ell} \in S$.
To do this we invoke \eqref{eq:aux3}. It suffices to show
\begin{equation}         \label{eq:aux6a}
  k-m \leq |i-\ell|.
\end{equation}
Observe $i< \ell$ by \eqref{eq:delta2aux1} and since $m+1 \leq k$; 
therefore $|i-\ell|=\ell-i$. 
By this and \eqref{eq:delta2aux1} we obtain \eqref{eq:aux6a},
and \eqref{eq:sp1b} follows.
Next we show \eqref{eq:delta2}.
By construction $k>m\geq r$ so $|r-k|=k-r$.
First assume 
\begin{equation}     \label{eq:delta2aux2}
   \ell+m < k.
\end{equation}
In this case $\ell+m-r \leq d$ so 
$|\Delta| = d-\ell-m+r$ and \eqref{eq:delta2} follows.
Next assume
\[
  k \leq \ell+m < d+r.
\]
In this case $|\Delta| = d-k+r+1$ and \eqref{eq:delta2} follows.
Next assume
\[
      d+r \leq \ell+m.
\]
In this case $|\Delta| = \ell-k+m+1$
and \eqref{eq:delta2} follows.
In any case \eqref{eq:delta2} holds and \eqref{eq:Aell} follows. 
This completes the proof.
\end{proof}

\begin{corollary}     \label{cor:EsrDDsDDss}
With reference to Definition {\rm \ref{def:T}},
for $0 \leq r,s \leq d$ we have
\[
  e^*_r D D^* D e^*_s =
    \sum_{t=0}^{\lfloor (r+s)/2 \rfloor} 
        e^*_r D e^*_t D e^*_s.
\]
\end{corollary}

\begin{proof}
Apply $\pi$ to  $D \otimes D^* \otimes D$, and
evaluate the result using Definition \ref{def:pi}, 
Lemma \ref{lem:kernel}, and Proposition \ref{prop:EsrDDsDEss}.
\end{proof}

\begin{proposition}    \label{prop:main}    \samepage
With reference to Notation {\rm \ref{notation}},
for an integer $n \geq 0$
the space
\begin{equation}   \label{eq:Es0DDsDEs0}
  e^*_0 D D^* D D^* D \cdots D D^* D e^*_0
    \qquad\qquad (\text{$n$ $D$'s})
\end{equation}
is equal to
\begin{equation}   \label{eq:Es0DEs0DEs0}
  e^*_0 D e^*_0 D e^*_0 D \cdots D e^*_0 D e^*_0
        \qquad\qquad (\text{$n$ $D$'s}).
\end{equation}
\end{proposition}

\begin{proof}
We assume $n \geq 1$; otherwise there is nothing to prove.
Since $\{e^*_i\}_{i=0}^d$ is a basis for $D^*$
it suffices to show that
\begin{equation}            \label{eq:Es0DEst1DEs0}
   e^*_0 D e^*_{t_1} D e^*_{t_2} D \cdots 
             D e^*_{t_n} D e^*_0
\end{equation}
is contained in \eqref{eq:Es0DEs0DEs0} for all sequences 
$\{t_i\}_{i=1}^n$ such that $0 \leq t_i  \leq d$ $(1 \leq i \leq n)$.
For each such sequence $\{t_i\}_{i=1}^n$ we define the {\em weight} to
be $\sum_{i=1}^n t_i$.
Suppose there exists a sequence $\{t_i\}_{i=1}^n$
such that \eqref{eq:Es0DEst1DEs0} is not contained in \eqref{eq:Es0DEs0DEs0}.
Of all such sequences,
pick one with minimal weight. Denote this weight by $w$.
For notational convenience define $t_0=0$ and $t_{n+1}=0$. 
We cannot have
$t_{i-1} \leq t_i$ for all $i$ $(1 \leq i \leq n+1)$;
otherwise $t_i=0$ for $0 \leq i \leq n+1$
and \eqref{eq:Es0DEst1DEs0} equals \eqref{eq:Es0DEs0DEs0}. 
By this and since $t_0=0$, there exists
$i$ $(1 \leq i \leq n)$ such that $t_{i-1} \leq t_i > t_{i+1}$.
Abbreviate $m = \lfloor (t_{i-1}+t_{i+1})/2 \rfloor$ and note
that $m< t_i$.
By Corollary \ref{cor:EsrDDsDDss} the space \eqref{eq:Es0DEst1DEs0} is 
contained in the space
\begin{equation}          \label{eq:sumEs0DEst1EsDEs0}
 \sum_{k=0}^m
  e^*_0 D e^*_{t_1} D \cdots D 
  e^*_{t_{i-1}} D e^*_k D e^*_{t_{i+1}} \cdots D e^*_{t_n} D e^*_0.
\end{equation}
For $0 \leq k \leq m$ the $k$-summand in \eqref{eq:sumEs0DEst1EsDEs0}
has weight less than $w$,
so this summand is contained in \eqref{eq:Es0DEs0DEs0}. 
Therefore \eqref{eq:Es0DEst1DEs0} is contained in \eqref{eq:Es0DEs0DEs0}, 
for a contradiction. The result follows.
\end{proof}

\begin{proofof}{Theorem \ref{thm:es0Tes0}}
(i):
Immediate from Proposition \ref{prop:main}
and since $T$ is generated by $D,D^*$.

(ii):
For $X,Y \in D$ we show $e^*_0Xe^*_0$, $e^*_0Ye^*_0$  commute.  
Consider the map $\pi$ from Definition \ref{def:pi},
and the space $R$ from Definition \ref{def:R}, where we take $r=0$, $s=0$.
By \cite[Theorem 8.3(ii)]{NT:eddde} and Lemma \ref{lem:injection}
we find
\[
  X \otimes e^*_0 \otimes Y - Y \otimes e^*_0 \otimes X \in R.
\]
In the above line we apply the map $\pi$ and use lemma \ref{lem:kernel}
to find
\[
  e^*_0 X e^*_0 Y e^*_0 = e^*_0 Y  e^*_0 X e^*_0.
\]
By this and since $e^{*2}_0=e^*_0$ the elements
$e^*_0Xe^*_0$, $e^*_0Ye^*_0$  commute.  

(iii):
Immediate from (i), (ii) above.
\end{proofof}

\section{The proof of Theorems \ref{thm:sharp}--\ref{thm:determined}}

\indent
In this section we use Theorem \ref{thm:es0Tes0} to prove 
Theorems \ref{thm:sharp}--\ref{thm:determined}.

\medskip

\begin{proofof}{Theorem \ref{thm:sharp}}
With reference to Notation \ref{notation},
assume $\mathbb{K}$ is algebraically closed, and
let $\Phi=(A;\{E^*_i\}_{i=0}^d;A^*$; $\{E^*_i\}_{i=0}^d)$ denote 
a tridiagonal system on $V$.
Consider the $T$-module structure on $V$ from Lemma \ref{lem:Tmodule}.
By construction $e^*_0V$ is nonzero, finite-dimensional, 
and invariant under $e^*_0Te^*_0$. 
Since $e^*_0Te^*_0$ is commutative by Theorem \ref{thm:es0Tes0}(iii),
and $\mathbb{K}$ is algebraically closed, 
there exists a nonzero $v \in e^*_0V$ which is a common eigenvector
for $e^*_0Te^*_0$. 
The subspace $Tv$ is nonzero and $T$-invariant so $Tv=V$
by the irreducibility of the $T$-module $V$.
In the equation $Tv=V$ we apply $e^*_0$ to both sides and use
$v=e^*_0v$ to get $e^*_0Te^*_0V=e^*_0V$.
Since $v$ is a common eigenvector for $e^*_0Te^*_0$,
we have
$e^*_0Te^*_0v=\text{Span}\{v\}$. By these comments the dimension
of $e^*_0V$ is 1. By construction $e^*_0V=E^*_0V$ so
$\rho_0=1$. The result follows.
\end{proofof}

\medskip

Our next goal is to prove Theorem \ref{thm:determined}. 
We start with a few definitions.

\medskip

\begin{definition}   \label{def:isoTDpair}   \samepage
Let $A,A^*$ denote a tridiagonal pair on $V$ and let
$A',A^{*\prime}$ denote a tridiagonal pair on $V'$.
We say these tridiagonal pairs are {\em isomorphic}
whenever there exists an isomorphism of $\mathbb{K}$-vector spaces
$\gamma : V \to V'$ such that
$\gamma A = A' \gamma$ and $\gamma A^* = A^{*\prime} \gamma$.
\end{definition}

\begin{definition}   \label{def:isoTDsystem}   \samepage
Let $\Phi=(A;\{E_i\}_{i=0}^d;A^*;\{E^*_i\}_{i=0}^d)$ denote a
tridiagonal system on $V$ and let
$\Phi'=(A';\{E^{\prime}_i\}_{i=0}^d; A^{*\prime};\{E^{*\prime}_i\}_{i=0}^d)$ 
denote a tridiagonal system on $V'$.
We say $\Phi$ and $\Phi'$ are {\em isomorphic}
whenever there exists an isomorphism of $\mathbb{K}$-vector spaces
$\gamma : V \to V'$ such that
$\gamma A = A' \gamma$, $\gamma A^* = A^{*\prime} \gamma$
and 
$\gamma E_i = E'_i \gamma$, $\gamma E^*_i=E^{*\prime}_i\gamma$ for
$0 \leq i \leq d$.
\end{definition}

\medskip

We now recall the split decomposition and the split sequence.
With reference to Notation \ref{notation},
for $0 \leq i \leq d$ we define
\begin{equation}          \label{eq:defUi}
   U_i = (E^*_0V+E^*_1V+\cdots+E^*_iV) \cap(E_iV+E_{i+1}V+\cdots+E_dV).
\end{equation}
By \cite[Theorem 4.6]{ITT}
\[
   V = U_0+U_1+\cdots+U_d  \qquad\qquad (\text{direct sum}),
\]
and for $0 \leq i \leq d$ both
\begin{align}
  U_0+U_1+\cdots+U_i &= E^*_0V+E^*_1V+\cdots+E^*_iV,    \label{eq:sumU0Ui} \\
  U_i+U_{i+1}+\cdots+U_d &= E_iV+E_{i+1}V+\cdots+E_dV.  \label{eq:sumUiUd}
\end{align}
By \cite[Corollary 5.7]{ITT}  $U_i$ has dimension $\rho_i$
for $0 \leq i \leq d$, where $\{\rho_i\}_{i=0}^d$ is the shape of $\Phi$.
By \cite[Theorem 4.6]{ITT} both
\begin{align}
  (A-\th_i I)U_i & \subseteq U_{i+1},    \label{eq:up}  \\
  (A^*-\th^*_i I)U_i & \subseteq U_{i-1}  \label{eq:down}
\end{align}
for $0 \leq i \leq d$, where $U_{-1}=0$ and $U_{d+1}=0$.
The sequence $\{U_i\}_{i=0}^d$ is called the {\em $\Phi$-split decomposition}
of $V$ \cite[Section 4]{ITT}.
Now assume $\Phi$ is sharp, so that $U_0$ has dimension $1$.
For $0 \leq i \leq d$ the space $U_0$ is invariant under
\begin{equation}        \label{eq:updown}
 (A^*-\th^*_1I)(A^*-\th^*_2I)\cdots(A^*-\th^*_iI)
 (A-\th_{i-1}I)\cdots(A-\th_1I)(A-\th_0I);
\end{equation}
let $\zeta_i$ denote the corresponding eigenvalue.
Note that $\zeta_0=1$.
We call the sequence $\{\zeta_i\}_{i=0}^d$ the {\em split sequence}
of $\Phi$.

\medskip

\begin{definition} {\rm \cite[Definition 14.5]{IT:Krawt}}
                         \label{def:parray}   \samepage
With reference to Notation \ref{notation}
assume $\Phi$ is sharp.
By the {\em parameter array} of $\Phi$ we mean the sequence
  $(\{\th_i\}_{i=0}^d; \{\th^*_i\}_{i=0}^d; \{\zeta_i\}_{i=0}^d)$
where $\{\zeta_i\}_{i=0}^d$ is the split sequence of $\Phi$.
We remark that the parameter array of $\Phi$ is defined only 
when $\Phi$ is sharp.
\end{definition}

\begin{lemma}      \label{lem:zetai}   \samepage
With reference to Notation {\rm \ref{notation}}, assume $\Phi$ is sharp
and let $\{\zeta_i\}_{i=0}^d$ denote the split sequence of $\Phi$.
Then for $0 \leq i \leq d$
the eigenvalue of 
\[
   E^*_0(A-\th_{i-1}I) \cdots (A-\th_1I)(A-\th_0I) E^*_0
\]
on the eigenspace $E^*_0V$ is   
\[
   \frac{\zeta_i}
        {(\th^*_0-\th^*_1)(\th^*_0-\th^*_2)\cdots(\th^*_0-\th^*_i)}.
\]
\end{lemma}

\begin{proof}
Let $\{U_i\}_{i=0}^d$ denote the $\Phi$-split decomposition of $V$.
Note that $U_0=E^*_0V$ so 
\begin{equation}    \label{eq:zetaiaux}
 (A^*-\th^*_1I)(A^*-\th^*_2I)\cdots(A^*-\th^*_iI) 
   (A-\th_{i-1}I) \cdots (A-\th_1I)(A-\th_0I) = \zeta_iI
\end{equation}
on $E^*_0V$.
To get the result, multiply both sides of \eqref{eq:zetaiaux}
on the left by $E^*_0$ and use $E^*_0A^*= \th^*_0E^*_0$.
\end{proof}

\begin{proofof}{Theorem \ref{thm:determined}}
Let
$\Phi=(A;\{E_i\}_{i=0}^d;A^*;\{E^*_i\}_{i=0}^d)$
and
$\Phi'=(A',\{E^{\prime}_i\}_{i=0}^d;A^{*\prime}$; $\{E^{*\prime}_i\}_{i=0}^d)$
denote sharp tridiagonal systems over $\mathbb{K}$ 
with the same parameter array 
$(\{\th_i\}_{i=0}^d$; $\{\th^*_i\}_{i=0}^d; \{\zeta_i\}_{i=0}^d)$.
We show that $\Phi$ and $\Phi'$ are isomorphic.
To this end, we consider the $\mathbb{K}$-algebra $T$ from
Definition \ref{def:T}.
Let $V$ (resp. $V'$) denote the vector space underlying $\Phi$ (resp. $\Phi'$).
Consider the $T$-module structure on $V$ and $V'$ from Lemma \ref{lem:Tmodule}.
It suffices to show that these $T$-modules are isomorphic.
Since $\Phi$ and $\Phi'$ are sharp, $e^*_0V$ and $e^*_0V'$ have dimension $1$; 
pick nonzero vectors $v \in e^*_0V$ and $v' \in e^*_0V'$. 
Define
$J = \{ x \in T \,|\, xv=0\}$,
$J' = \{ x \in T \,|\, xv'=0\}$
and note that $J$, $J'$ are left ideals of $T$.
The quotient vector space $T/J$ has a natural $T$-module structure and the map
\begin{align*}
   T/J  \quad & \to  \quad  V   \\
  x + J \quad & \mapsto   \quad xv
\end{align*}
is an isomorphism of $T$-modules. A similar result holds for $J'$. 
Therefore, is suffices to show that $J=J'$.
By Theorem \ref{thm:es0Tes0}, Lemma \ref{lem:zetai}, 
and since $\Phi$, $\Phi'$ have the same parameter array,
\begin{equation}                  \label{eq:Jcapes0Tes0}
  J \cap e^*_0Te^*_0 = J' \cap e^*_0Te^*_0.
\end{equation}
We claim  $J \subseteq J'$.
To prove the claim, consider the subspace $Jv'$ of $V'$.
By construction $TJ \subseteq J$ so $Jv'$ is a $T$-submodule of $V'$. 
Therefore $Jv'=0$ or $Jv'=V'$. 
Assume $Jv'=V'$; otherwise $J \subseteq J'$ and the claim is proved.
Since $Jv'=V'$ there exists $x \in J$ such that $xv'=v'$.
By construction $e^*_0 x e^*_0 \in e^*_0Te^*_0$.
We have 
$e^*_0 x e^*_0 v  
             = e^*_0 x v 
             = e^*_0 0
             = 0$
so $e^*_0 x e^*_0 \in J$.
Also $e^*_0 x e^*_0 v'  
             = e^*_0 x v' 
             = e^*_0 v' 
             = v'$
so $e^*_0 x e^*_0 \not\in J'$.
This contradicts \eqref{eq:Jcapes0Tes0} so we must have $J \subseteq J'$.
By symmetry we have $J' \subseteq J$ so $J=J'$ and
the theorem follows.
\end{proofof}

\medskip

We mention a few variations on Theorem \ref{thm:determined}.

\medskip

\begin{corollary}   \label{cor:determined}
Assume $\mathbb{K}$ is algebraically closed. 
Then two tridiagonal systems over $\mathbb{K}$ are isomorphic
if and only if they have the same parameter array.
\end{corollary}

\begin{proof}
Combine Theorem \ref{thm:sharp} and Theorem \ref{thm:determined}.
\end{proof}

\begin{definition}   \samepage
Let $A,A^*$ denote a tridiagonal pair.
By a {\em parameter array} of $A,A^*$ we mean the parameter
array of an associated tridiagonal system.
\end{definition}

\begin{theorem}    \label{thm:determined2}     \samepage
Two sharp tridiagonal pairs over $\mathbb{K}$ are isomorphic if and only if they
have a parameter array in common.
\end{theorem}

\begin{proof}
Let $A,A^*$ and $A', A^{*\prime}$ denote sharp tridiagonal pairs over 
$\mathbb{K}$ that have a parameter array in common. By construction this
is the parameter array for a tridiagonal system
$(A;\{E_i\}_{i=0}^d; A^*; \{E^*_i\}_{i=0}^d)$ 
and a tridiagonal system
$(A';\{E'_i\}_{i=0}^d; A^{*\prime};\{E^{*\prime}_i\}_{i=0}^d)$.
By Theorem \ref{thm:determined} these tridiagonal systems are isomorphic
so the tridiagonal pairs $A,A^*$ and $A',A^{*\prime}$ are isomorphic.
This proves one direction of the theorem; the other
direction is clear.
\end{proof}

\begin{corollary}   \samepage
Assume $\mathbb{K}$ is algebraically closed.
Then two tridiagonal pairs over $\mathbb{K}$ are isomorphic
if and only if they have a parameter array in common.
\end{corollary}

\begin{proof}
Combine Theorem \ref{thm:sharp} and Theorem \ref{thm:determined2}.
\end{proof}

\medskip

We now turn our attention to Theorems \ref{thm:bilin} and
\ref{thm:existanti}.

\medskip

A map $\b{\;,\,} : V \times V \to \mathbb{K}$ is called a
{\em bilinear form} on $V$ whenever the following conditions hold
for $u,v,w \in V$ and for $\alpha \in \mathbb{K}$:
(i) $\b{u+v,w}=\b{u,w}+\b{v,w}$;
(ii) $\b{\alpha u,v}=\alpha \b{u,v}$;
(iii) $\b{u,v+w}=\b{u,v}+\b{u,w}$;
(iv) $\b{u,\alpha v}=\alpha \b{u,v}$.
We observe that a scalar multiple of a bilinear form is a bilinear form.
Let $\b{\;,\,}$ denote a bilinear form on $V$.
Then the following are equivalent:
(i) there exists a nonzero $u \in V$ such that $\b{u,v}=0$ for all $v\in V$;
(ii) there exists a nonzero $v \in V$ such that $\b{u,v}=0$ for all $u\in V$.
The form $\b{\;,\,}$ is said to be {\em degenerate} whenever (i), (ii) hold
and {\em nondegenerate} otherwise.
The form is said to be {\em symmetric} whenever $\b{u,v}=\b{v,u}$
for all $u,v \in V$.

\medskip

\begin{proofof}{Theorem \ref{thm:bilin}}
Follows from Theorem \ref{thm:determined} and
\cite[Theorem 11.4]{NT:sharp}.
\end{proofof}

\medskip

By an {\em anti-automorphism} of $\text{End}(V)$ we mean an isomorphism
of $\mathbb{K}$-vector spaces $\sigma: \text{End}(V) \to \text{End}(V)$
such that $(XY)^\sigma=Y^\sigma X^\sigma$ for all $X,Y \in \text{End}(V)$.

\medskip

\begin{proofof}{Theorem \ref{thm:existanti}}
Follows from Theorem \ref{thm:determined} and
\cite[Theorem 11.5]{NT:sharp}.
\end{proofof}

\medskip

We mention a few variations on Theorems \ref{thm:bilin}, \ref{thm:existanti}.

\medskip

\begin{corollary}    \samepage
Assume $\mathbb{K}$ is algebraically closed and let $A,A^*$ denote a
tridiagonal pair on $V$.
Then there exists a nonzero bilinear form $\b{\;,\,}$ on $V$ such that
$\b{Au,v}=\b{u,Av}$ and $\b{A^*u,v}=\b{u,A^*v}$ for all $u,v \in V$.
This form is unique up to multiplication by a nonzero scalar in $\mathbb{K}$.
This form is nondegenerate and symmetric.
\end{corollary}

\begin{proof}
Combine Theorem \ref{thm:sharp} and Theorem \ref{thm:bilin}.
\end{proof}

\begin{corollary}   \samepage
Assume $\mathbb{K}$ is algebraically closed and let $A,A^*$ denote a
tridiagonal pair on $V$.
Then there exists a unique anti-automorphism $\dagger$ of $\text{\rm End}(V)$
that fixes each of $A,A^*$.
Moreover $X^{\dagger\dagger}=X$ for all $X \in \text{\rm End}(V)$.
\end{corollary}

\begin{proof}
Combine Theorem \ref{thm:sharp} and Theorem \ref{thm:existanti}.
\end{proof}

\bigskip

{\small

\bibliographystyle{plain}

}

\bigskip\bigskip\noindent
Kazumasa Nomura\\
College of Liberal Arts and Sciences\\
Tokyo Medical and Dental University\\
Kohnodai, Ichikawa, 272-0827 Japan\\
email: knomura@pop11.odn.ne.jp

\bigskip\noindent
Paul Terwilliger\\
Department of Mathematics\\
University of Wisconsin\\
480 Lincoln Drive\\ 
Madison, Wisconsin, 53706 USA\\
email: terwilli@math.wisc.edu

\bigskip\noindent
{\bf Keywords.}
Leonard pair, tridiagonal pair, $q$-Racah polynomial, orthogonal polynomial.

\noindent
{\bf 2000 Mathematics Subject Classification}.
05E35, 05E30, 33C45, 33D45.

\end{document}